\documentclass[own,a4paper]{mysvjour}

\newif\ifcenterpage
\centerpagetrue

\ifcenterpage
\makeatletter
\@settopoint\oddsidemargin \@settopoint\marginparwidth
\@settopoint\evensidemargin
\makeatother
\fi

\newif\ifusepstricks   

\usepstrickstrue   

\ifusepstricks
  \usepackage{pstricks,pst-node,pstcol,pst-node,pst-eps} 
\fi

\usepackage{amssymb}
\usepackage{amsmath}
\usepackage{bbm}
\usepackage{booktabs}
\usepackage{enumerate}
\usepackage{epsfig}
\usepackage{subfigure}
\usepackage{times}
\usepackage{url}

\newcommand{\NP}{\ensuremath{\mathcal{NP}}}
\newcommand{\coNP}{\ensuremath{\text{\rm co-}\mathcal{NP}}}

\newcommand{\R}{\ensuremath{\mathbb{R}}}

\DeclareMathOperator\convOp{conv}

\newcommand{\VFI}[1]{\ensuremath{\mathcal{I}({#1})}}
\newcommand{\aff}[1]{\ensuremath{\operatorname{aff}({#1})}}
\newcommand{\bary}[1]{\ensuremath{\beta({#1})}}

\newcommand{\coordVertices}[1]{\ensuremath{\operatorname{vert}(#1)}}
\newcommand{\dualPolytope}[1]{\ensuremath{#1^{\star}}}
\newcommand{\edgebij}[2]{\ensuremath{\Psi_{{#1},{#2}}}}
\newcommand{\edges}[1]{\ensuremath{\operatorname{E}(#1)}}
\newcommand{\graph}[1]{\ensuremath{\mathcal{G}({#1})}}

\newcommand{\neighborsIncl}[1]{\ensuremath{\overline{\operatorname{N}}(#1)}}
\newcommand{\neighbors}[1]{\ensuremath{\text{N}(#1)}}
\newcommand{\order}[1]{\ensuremath{\mathcal{O}(#1)}}
\newcommand{\poly}[1]{\ensuremath{\mathcal{P}({#1})}}

\newcommand{\vertices}[1]{\ensuremath{\operatorname{V}(#1)}}

\newcommand{\polymake}{\texttt{polymake}}
\newcommand{\nauty}{\texttt{nauty}}
\newcommand{\javaview}{\texttt{javaview}}

\newcommand{\topic}[1]{\subsubsection*{#1.} }

\journalname{Graphs and Combinatorics}
\begin{document}

\title{On the Complexity of\\Polytope Isomorphism Problems}
\author{Volker Kaibel%
           \thanks{Supported by a DFG 
                   Gerhard-Hess-Forschungsf\"orderungspreis donated to
                   G\"unter M.~Ziegler (Zi 475/2-3).}
        \and Alexander Schwartz}

\institute{%
    Technische Universit\"at Berlin\\
    Fakult\"at~II, Institut f\"ur Mathematik, MA 6--2\\
    Stra\ss e des 17. Juni~136\\10623~Berlin, Germany\\
    \email{\{schwartz,kaibel\}@math.tu-berlin.de}%
}

\date{July 2002\\[.5cm]\emph{(To appear in Graphs and Combinatorics.)}}
\maketitle


\begin{abstract}
  We show that the problem to decide whether two (convex) polytopes,
  given by their vertex-facet incidences, are combinatorially
  isomorphic is graph isomorphism complete, even for simple or
  simplicial polytopes.  On the other hand, we give a polynomial time
  algorithm for the combinatorial polytope isomorphism problem in
  bounded dimensions.  Furthermore, we derive that the problems to
  decide whether two polytopes, given either by vertex or by facet
  descriptions, are projectively or affinely isomorphic are graph
  isomorphism hard.
\end{abstract}




\medskip
{\small
\noindent\textsf{\textbf{MSC~2000:}}
  52B05, 05C60, 52B11, 68R10}



\section{{Introduction}}

\topic{Combinatorial isomorphism of polytopes}
When treated as a combinatorial object, a polytope (i.e., a bounded
convex polyhedron) is identified with its \emph{face lattice}, i.e.,
the lattice formed by its faces, which are ordered by inclusion. Two
polytopes are considered \emph{combinatorially isomorphic} if their
face lattices are isomorphic, i.e., if there is an in both directions
inclusion preserving bijection between their sets of faces.

Since the face lattice is both atomic and coatomic, the entire
combinatorial structure of a polytope~$P$ is encoded in its
\emph{vertex-facet incidences}, i.e., in a bipartite graph~$\VFI{P}$,
whose two shores represent the vertices and facets, where an edge
indicates that the vertex corresponding to the one end node is
contained in the facet corresponding to the other one. The
\emph{polytope isomorphism problem} is the problem to decide whether
two polytopes~$P$ and~$Q$, given by their vertex-facet incidences, are
combinatorially isomorphic.  For all concepts and notations concerning
polytope theory, we refer to Ziegler's book~\cite{Ziegler95}.

Two graphs~$G$ and~$G'$ are called \emph{isomorphic} if there is a
bijection~$\varphi$ between their vertex sets such that
$\{\varphi(v),\varphi(w)\}$ is an edge of~$G'$ if and only if
$\{v,w\}$ is an edge of~$G$. The problem to decide whether two
graphs~$G$ and~$G'$ are isomorphic is known as the \emph{graph
  isomorphism problem}.  If the nodes and edges of the two graphs are
labeled, and if we consider only isomorphisms that preserve the
labels, then the corresponding decision problem is called the
\emph{labeled graph isomorphism problem}. Clearly, the (usual) graph
isomorphism problem is the special case of the labeled graph
isomorphism problem, where all nodes have the same label. Conversely,
the labeled graph isomorphism problem can be Karp-reduced (see below)
to the graph isomorphism problem.

The polytope isomorphism problem for two polytopes~$P$ and~$Q$ is a
special case of the labeled graph isomorphism problem, where in the
bipartite graphs~$\VFI{P}$ and~$\VFI{Q}$ each node is either
labeled ``vertex'' or ``facet.''

Since every isomorphism between two polytopes or between two graphs is
determined by its restriction to the vertices or nodes, respectively,
we will often identify an isomorphism with that restriction. Thus, two
polytopes are combinatorially isomorphic if there is a bijection
between their sets of vertices that induces a bijection between their
sets of facets.

\topic{Checking isomorphism by computer} 
For computer systems dealing with (the combinatorial structures of)
polytopes, the problem of checking two polytopes for combinatorial
isomorphism is quite important. Gr\"unbaum writes in his classical
book \cite[footnote on p.~39]{Gruenbaum67}: 
\begin{quote}
  \textit{
    For two given polytopes it is, in principle, easy to determine
    whether they are combinatorially equivalent or not. It is enough
    (\dots) to check whether there exists any inclusion preserving
    one-to-one correspondence between the two sets of faces. However,
    this procedure is practically feasible only if the number of faces
    is rather small.
  }
\end{quote}

The \polymake -system of Gawrilow and Joswig~\cite{GawrilowJoswig00b}
currently implements the isomorphism test for two polytopes~$P$
and~$Q$ by checking whether the bipartite graphs $\VFI{P}$ and
$\VFI{Q}$ are isomorphic as labeled graphs. This is done by using the
software package \texttt{nauty} by McKay~\cite{McKay81}.

One of our results (Theorem~\ref{thm:polyIsoCompl}) shows that in
order to solve the general polytope isomorphism problem the only way
is indeed to use an algorithm for the general graph isomorphism
problem. Our second main result (Theorem~\ref{thm:iso of d-polytopes})
shows that one might take advantage of the fact that the polytopes
dealt with in computer systems usually have rather small dimensions.

\topic{Geometric notions of isomorphism of polytopes}
If one is  concerned with coordinate representations of
polytopes then geometric notions of isomorphism become important. In
particular, two polytopes $P\subset\R^{d'}$, $Q\subset\R^{d''}$ are
called \emph{affinely} or \emph{projectively} \emph{isomorphic} if
there is an affine or projective, respectively, map from~$\R^{d'}$ to
$\R^{d''}$ inducing a bijection between~$P$ to~$Q$.  Two polytopes are
\emph{congruent} if there is an affine isomorphism between them that
is induced by an orthogonal matrix.

The decision problems corresponding to these geometric notions of
isomorphism are the \emph{polytope congruence problem} and the
\emph{affine} or \emph{projective polytope isomorphism problem}, where
for each of them the two polytopes either are given by vertex
coordinates ($\mathcal{V}$-descriptions) or by inequality coefficients
($\mathcal{H}$-descriptions).  

The chain of implications ``congruent $\Rightarrow$ affinely
isomorphic $\Rightarrow$ projectively isomorphic $\Rightarrow$
combinatorially isomorphic'' is well-known.  We will mainly be
concerned with combinatorial isomorphism; nevertheless, our results
also have implications for geometric isomorphism (see
Theorem~\ref{thm:geomiso}).

\topic{Graph isomorphism and reductions} 
The complexity status of the general graph isomorphism problem is open.
While it is obvious that the problem is contained in the complexity
class~$\NP$, all attempts either to show that it is also contained
in~$\coNP$ (or even that it can be solved in polynomial time) as well
as all efforts in the direction of proving its $\NP$-completeness have
failed so far. In fact, this apparent difficulty of classifying the
complexity is shared by a number of isomorphism problems.

There are a variety of problems which are in a certain sense as
difficult as the graph isomorphism problem, which means that they are
efficiently reducible to the graph isomorphism problem and vice versa,
where two concepts of reducibility are important.

A decision problem~$A$ is \emph{Karp reducible} to another decision
problem $B$, if there is a polynomial time algorithm which constructs
from an instance~$I$ of~$A$ an instance $J$ of $B$ with the property
that the answer for $J$ is ``yes'' if and only if the answer for $I$
is ``yes.'' Two decision problems $A$ and $B$ are called \emph{Karp
  equivalent} if $A$ is Karp reducible to $B$ and vice versa.  A
(decision) problem which is Karp equivalent to the graph isomorphism
problem is called \emph{graph isomorphism complete}.  Often a decision
problem~$\Pi$ is called \emph{graph isomorphism hard} if the graph
isomorphism problem is Karp reducible to $\Pi$.

Among the graph isomorphism complete problems are the restriction of
the graph isomorphism problem to the class of bipartite graphs (and
therefore comparability graphs), regular graphs
\cite{Booth78,Miller79,CorneilKirkpatrick80}, line graphs (see Harary
\cite{Harary69}), chordal graphs \cite{BoothLueker75}, and
self-complementary graphs \cite{ColbournColbourn78}; even the question
to decide whether a graph is self-complementary is graph isomorphism
complete \cite{ColbournColbourn78}.  A recent hardness result on a
problem that is not a special case of the graph isomorphism problem is
due to Kutz \cite{Kutz02}: he proved that for every $k>0$, checking if
a directed graph has a ``$k$-th root'' is graph isomorphism hard.

Other graph isomorphism complete problems occur in algebra (semi-group
isomorphism \cite{Booth78}, finitely presented algebra isomorphism
\cite{Kozen77}), as well as in topology (homeomorphism of
2-complexes \cite{Shawe-TaylorPisanski94} and homotopy equivalence
\cite{ZalcsteinFranklin86}).
  
Some interesting problems related to graph isomorphism are equivalent
to graph isomorphism if we use a weaker concept of reducibility: A
problem~$A$ is called \emph{Turing reducible} to a problem $B$ if
there is a polynomial time algorithm for the problem $A$ that might
use an oracle for solving~$B$, where each call to the oracle is
assumed to take only one step.  Two problems~$A$ and~$B$ are
\emph{Turing equivalent} if $A$ is Turing reducible to $B$ and vice
versa.  Mathon \cite{Mathon79} proved that a number of problems on
graphs are Turing equivalent to graph isomorphism, including counting
the number of \emph{automorphisms} (i.e., isomorphisms between the
graph and itself), finding a set of generators of the automorphism
group, and constructing the \emph{automorphism partition} (i.e.~the
orbits of the nodes under the automorphism group).

Lubiw \cite{Lubiw81} showed that the problem to decide, whether for a
given graph~$G$ and two specified nodes~$v$ and~$w$ of~$G$ there is an
automorphism of~$G$ not mapping~$v$ to~$w$, is Turing equivalent to
the graph isomorphism problem.  In contrast to this, she showed that
deciding whether a graph has a fix-point free automorphism is
$\NP$-complete.

For further information about the graph isomorphism problem, we refer
to the books by Hoffmann \cite{Hoffmann82}, and by K\"obler,
Sch\"oning \&\ Tor\'{a}n~\cite{KoeblerSchoeningToran93}, as well
as to the surveys by Read \& Corneil \cite{ReadCorneil77}, Babai
\cite{Babai95}, and Fortin \cite{Fortin96}.

\topic{Overview} 
In~Sect.~\ref{sec:arbDim}, we show that the (combinatorial) polytope
isomorphism problem is graph isomorphism complete
(Theorem~\ref{thm:polyIsoCompl}). This remains true for \emph{simple}
(every vertex figure is a simplex) as well as for \emph{simplicial}
(every facet is a simplex) polytopes.  The general graph isomorphism
problem can be Karp-reduced to each of the three geometric polytope
isomorphism problems stated above (Theorem~\ref{thm:geomiso}).
Furthermore, the graph isomorphism problem restricted to graphs of
polytopes (formed by their vertices and edges) is graph isomorphism
complete, even for the graphs of simple polytopes and the graphs of
simplicial polytopes (Theorem~\ref{thm:polyGraphIsoCompl}).

In Sect.~\ref{sec:boundDim}, we describe a polynomial time algorithm
for the isomorphism problem of polytopes of bounded dimensions
(Theorem~\ref{thm:iso of d-polytopes}).



\section{Hardness Results for Arbitrary Dimension}
\label{sec:arbDim}

\begin{figure}[ht]
  \centering \hfill
    \subfigure[]{
      {\epsfig{figure=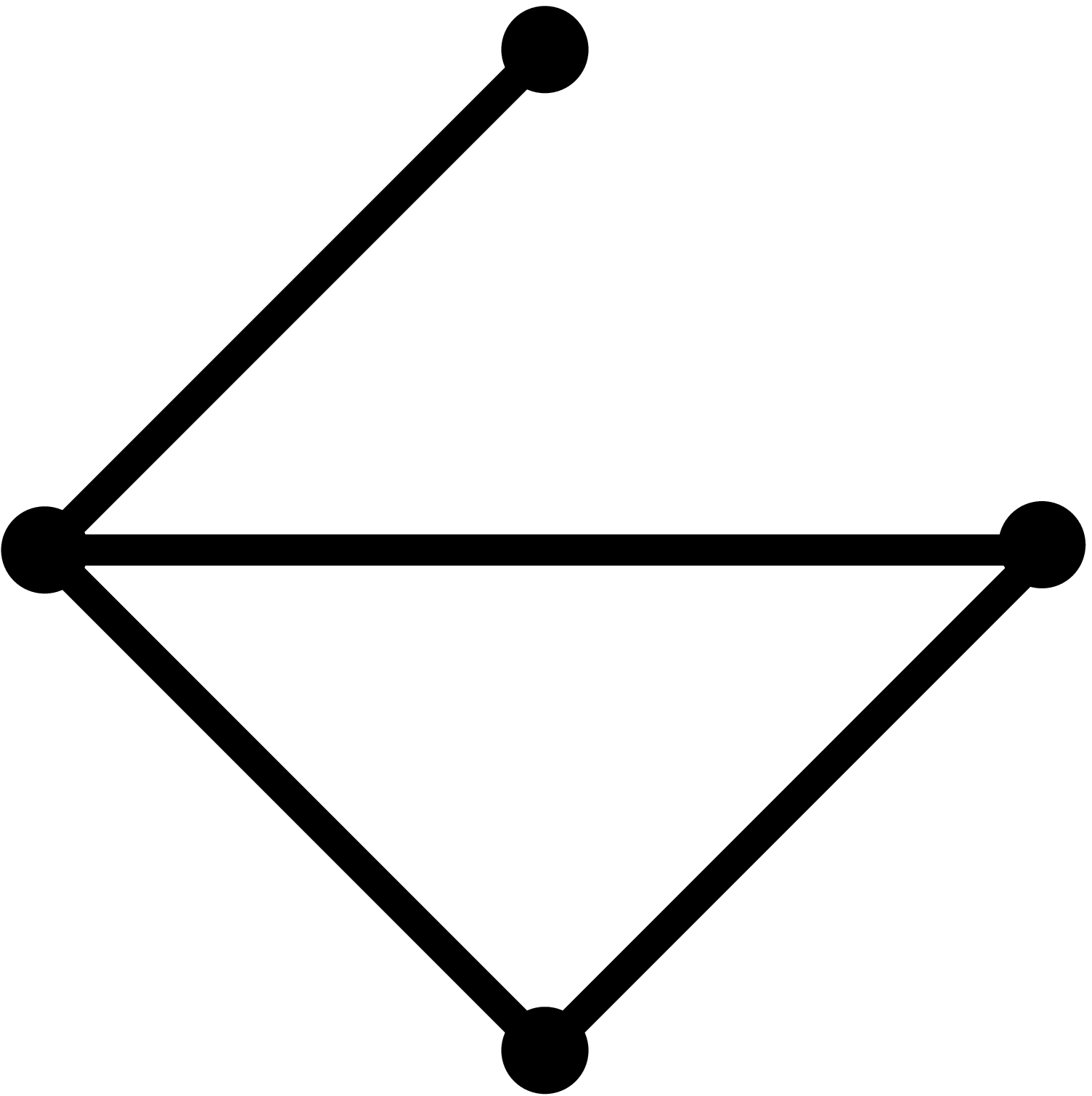,
        width=.14\textwidth}}}
    \hfill
    \subfigure[]{
      {\epsfig{figure=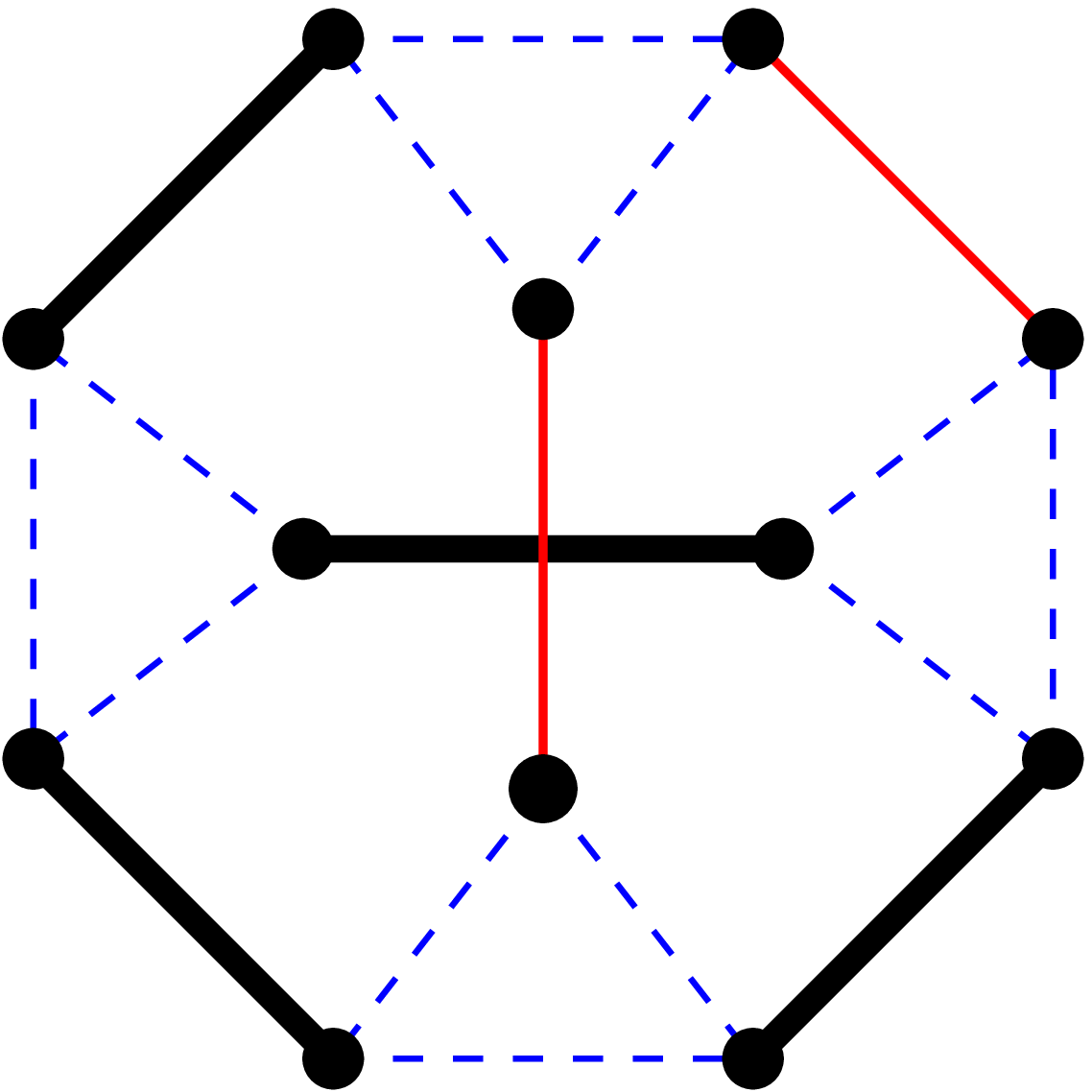,
        width=.2\textwidth}}}\hfill
    \subfigure[]{
      {\epsfig{figure=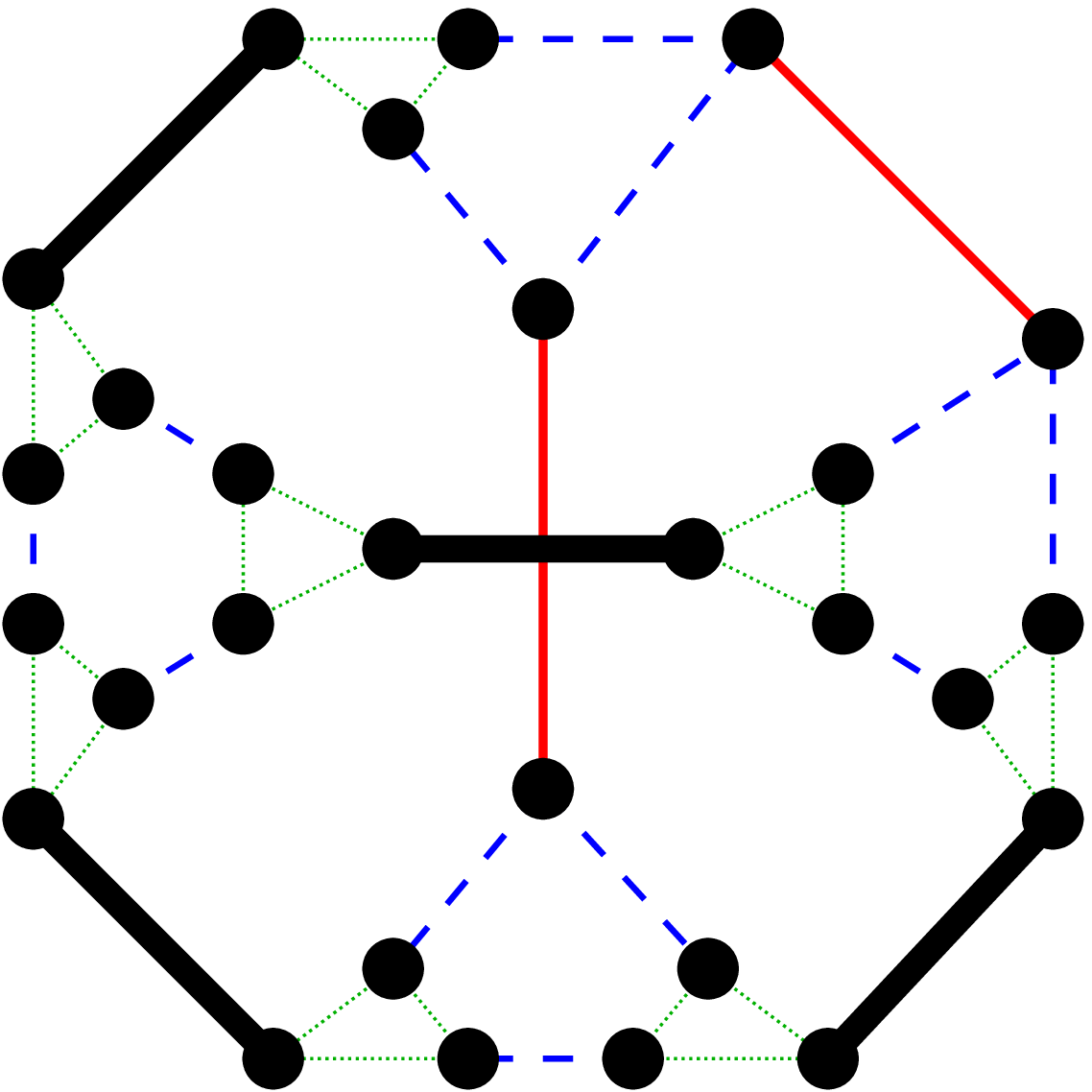,
        width=.2\textwidth}}}
    \hspace*{\fill}
    \subfigure{
      {\epsfig{figure=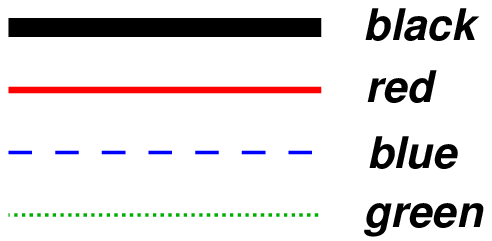,
        width=.065\textwidth,angle=90}}}
    \hspace*{\fill}
    \caption{Illustration of the construction of the
      graph~$\graph{\poly{G}}$ (depicted in~(c)) from a graph~$G$
      (depicted in~(a)). The graph of the intermediate
      polytope~$\Gamma_{n-1}$ is shown in~(b).}
    \label{fig:PfromG}
\end{figure}

\begin{figure}
  \centering
  \epsfig{figure=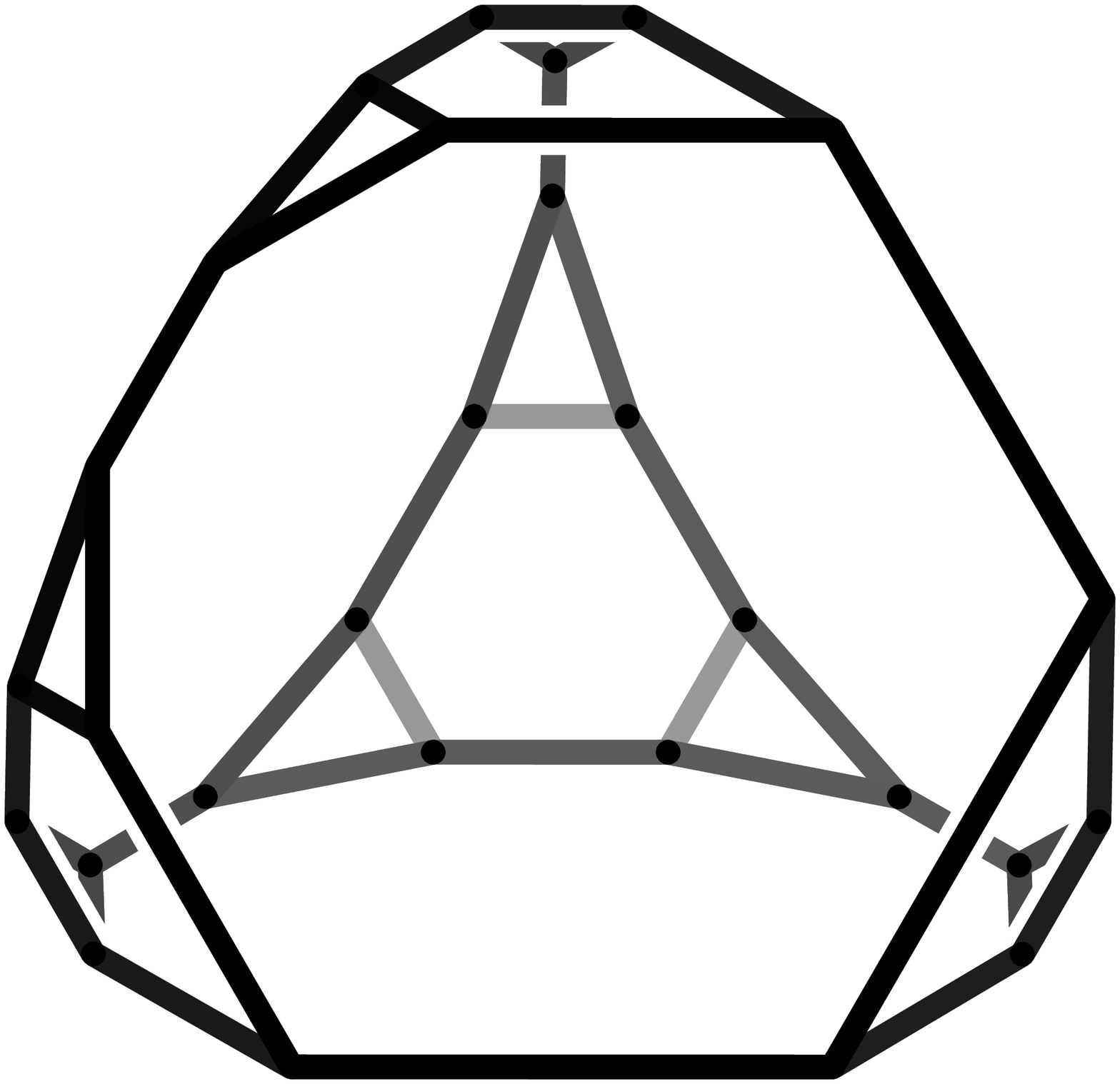,width=.3\textwidth}
  \caption{The polytope~$\poly{G}$ for the graph~$G$ in
    Fig.~\ref{fig:PfromG}.
    (Image produced with \polymake\ \cite{GawrilowJoswig00b} 
     and \javaview\ \cite{PolthierKhademPreussReitebuch01}.)}
  \label{fig:polyEx}
\end{figure}

The results in this section are based on the following construction
that produces from any given graph~$G=(V,E)$ on $n=|V|$ nodes a
certain polytope~$\poly{G}$ (see Fig.~\ref{fig:PfromG} and Fig.~\ref{fig:polyEx}). We denote
by~$\graph{P}$ the graph of a polytope~$P$ defined by the vertices and
one-dimensional faces of~$P$.

\begin{figure}
    \begin{center}
      \epsfig{figure=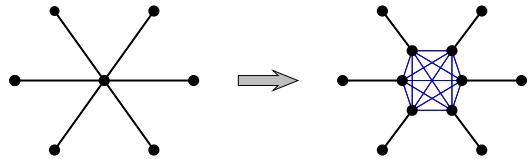,width=.6\textwidth}
      \caption{Cutting off a vertex~$v$ of a simple $k$-polytope
        means to replace the node~$v$ in the graph by a $k$-clique.}
      \label{fig:vertexClique}
    \end{center}
  \end{figure}

\begin{list}{}{}
  \item[\emph{First step.}] Choose an arbitrary bijection of~$V$ to the $n$
  vertices of an $(n-1)$-dimensional simplex~$\Delta_{n-1}$, thus
  embedding~$G$ into the graph~$\graph{\Delta_{n-1}}$
  of~$\Delta_{n-1}$. We call those edges of~$\graph{\Delta_{n-1}}$
  which are images under that embedding \emph{black} edges, and the
  other ones \emph{red} edges.
  
  \item[\emph{Second step.}] Cut off each vertex of the simplex~$\Delta_{n-1}$
  to obtain a polytope~$\Gamma_{n-1}$. The graph~$\graph{\Gamma_{n-1}}$ of~$\Gamma_{n-1}$
  arises from~$\graph{\Delta_{n-1}}$ by replacing each vertex by
  an~$(n-1)$-clique (see Fig.~\ref{fig:vertexClique}).  We call the
  edges of these cliques \emph{blue} edges. Thus, $\graph{\Gamma_{n-1}}$ has
  black edges corresponding to the edges of~$G$, red edges
  corresponding to the edges of the complement of~$G$, and blue edges
  coming from cutting off the vertices of~$\Delta_{n-1}$.

  \item[\emph{Third step.}] Construct~$\poly{G}$ from~$\Gamma_{n-1}$ by cutting off those 
  vertices that are incident to black edges. We call the edges of the 
  $(n-1)$-cliques that arise \emph{green} edges. 
\end{list}

The polytope~$\poly{G}$ is a simple $(n-1)$-polytope (see
Fig.~\ref{fig:polyEx}); its dual $\dualPolytope{\poly{G}}$ is a
simplicial $(n-1)$-polytope. The operation dual to cutting off a
vertex of a polytope~$P$ is a \emph{stellar subdivision} of the
corresponding facet of~$\dualPolytope{P}$, which means to place a new
vertex slightly beyond this facet. Thus, $\dualPolytope{\poly{G}}$ can
be obtained from a simplex by iteratively placing new vertices beyond
facets; it is a \emph{stacked polytope}. In the graph of a polytope,
placing a vertex beyond a facet has the effect that a new vertex is
added which is connected to all vertices of that facet.

\begin{remark}
\label{rem:GP}
  There are polynomial time algorithms that compute from a graph~$G$
  \begin{enumerate}[(i)]
  \item $\graph{\poly{G}}$ and $\graph{\dualPolytope{\poly{G}}}$,
  \item the vertex-facet incidences of~$\poly{G}$ and
    $\dualPolytope{\poly{G}}$,
  \item $\mathcal{V}$-descriptions of~$\poly{G}$ and
    $\dualPolytope{\poly{G}}$, and
  \item $\mathcal{H}$-descriptions  of~$\poly{G}$ and $\dualPolytope{\poly{G}}$.
  \end{enumerate}
\end{remark}

The important  property of~$\poly{G}$ is that it encodes the
entire structure of~$G$.

\begin{proposition}
\label{prop:GP}
For two graphs~$G$ and~$H$ on at least three nodes the
following five statements are equivalent.
\begin{enumerate}[(i)]
  \item $G$ is isomorphic to~$H$. 
  \item $\graph{\poly{G}}$ is isomorphic to~$\graph{\poly{H}}$. 
  \item $\graph{\dualPolytope{\poly{G}}}$ is isomorphic
        to~$\graph{\dualPolytope{\poly{H}}}$.
  \item $\poly{G}$ is isomorphic to $\poly{H}$.
  \item $\dualPolytope{\poly{G}}$ is isomorphic to
        $\dualPolytope{\poly{H}}$.
\end{enumerate}
\end{proposition}

\begin{proof}
  We start by proving the equivalence of~(i) and~(ii).

  Any isomorphism between two graphs~$G$ and~$H$ induces a color
  preserving isomorphism between the two complete graphs constructed
  from~$G$ and~$H$ in the first step. Of course, such a color
  preserving isomorphism induces a color preserving isomorphism of the
  graphs of the polytopes constructed in the second step, which
  finally gives rise to an isomorphism of the
  graphs~$\graph{\poly{G}}$ and~$\graph{\poly{H}}$ of the two
  polytopes constructed in the third step.
  
  In order to prove the converse direction, let~$G$ and~$H$ be two
  graphs on~$n$ and~$n'$ nodes ($n,n'\geq 3$), respectively, and
  let~$\varphi$ be an isomorphism between~$\graph{\poly{G}}$
  and~$\graph{\poly{H}}$.  Since~$\graph{\poly{G}}$ is
  $(n-1)$-regular and ~$\graph{\poly{H}}$ is $(n'-1)$-regular, we
  have $n=n'$. If $n=3$, then both~$\graph{\poly{G}}$
  and~$\graph{\poly{H}}$ are cycles of length~$\ell$. Since in this
  case, the number of edges of~$G$ as well as of~$H$ must be
  $(\ell-6)/2\in\{0,1,2,3\}$, $G$ and~$H$ are isomorphic. Thus, we
  may assume $n\geq 4$.
  
  We consider~$\graph{\poly{G}}$ and~$\graph{\poly{H}}$ colored as
  defined in the description of the construction. In both graphs, all
  $(n-1)$-cliques are node-disjoint. Each of these cliques either
  consists of green or of blue edges (blue cliques might arise from
  isolated nodes). Consider the graphs that arise
  from~$\graph{\poly{G}}$ and~$\graph{\poly{H}}$ by shrinking all
  $(n-1)$-cliques. Those nodes that come from shrinking green cliques
  are contained in (maximal) $(n-1)$-cliques in the shrunken graphs,
  while those coming from blue cliques are not (notice that for graphs
  without edges this statement indeed only holds for \emph{maximal}
  $(n-1)$-cliques). This shows that~$\varphi$ preserves the colors of
  $(n-1)$-cliques.
 
  Let~$G'$ and~$H'$ be the graphs that are obtained from shrinking the
  \emph{green} cliques in~$\graph{\poly{G}}$ and~$\graph{\poly{H}}$,
  respectively. Since~$\varphi$ maps green cliques to green cliques, it
  induces an isomorphism~$\psi$ between~$G'$ and~$H'$.  Since the
  shrinking operations do not generate multiple edges, the graphs~$G'$
  and~$H'$ inherit colorings of their edges from~$\graph{\poly{G}}$
  and~$\graph{\poly{H}}$, respectively. Because an edge
  of~$\graph{\poly{G}}$ or~$\graph{\poly{H}}$ is red if and only if
  it is not adjacent to a green edge, the isomorphism~$\psi$ preserves
  red edges.
  
  In the graphs~$G'$ and~$H'$ the only $(n-1)$-cliques are the ones
  formed by the blue edges. Again, these cliques are pairwise
  node-disjoint. Thus the isomorphism~$\psi$ between~$G'$ and~$H'$
  induces a color preserving isomorphism between the (complete) graphs
  obtained by shrinking all $(n-1)$-cliques in~$G'$ and~$H'$ (which,
  again, does not produce multiple edges). This, finally, yields
  an isomorphism between~$G$ and~$H$.
  
  The equivalence of~(ii) and~(iv) follows from a theorem of
  Blind and Mani~\cite{BlindMani87} (see also Kalai's beautiful
  proof~\cite{Kalai88}) stating that two simple polytopes are
  isomorphic if and only if their graphs are isomorphic. For the
  special polytopes arising from our construction, the equivalence
  can, however, be alternatively deduced similarly to the proof
  of ``(ii) $\Rightarrow$ (i).''

  Statements~(iv) and~(v) obviously are equivalent.
  
  Unlike the situation for simple polytopes, it is, in general, not
  true that two (simplicial) polytopes are isomorphic if and only
  their graphs are isomorphic. Nevertheless, for stacked polytopes
  like $\dualPolytope{\poly{G}}$ and $\dualPolytope{\poly{H}}$ it is
  true (this follows, e.g., from the fact that one can reconstruct the
  vertex-facet incidences of a stacked $d$-polytope from its graph by
  iteratively removing vertices of degree~$d$). Thus, finally the
  equivalence of~(iii) and~(v) is established.\hspace*{\fill}\qed
\end{proof}

The equivalences ``(i) $\Leftrightarrow$ (ii)'' and ``(i)
$\Leftrightarrow$ (iii)'' in Proposition~\ref{prop:GP} together with
part~(i) of Remark~\ref{rem:GP} immediately imply that the restriction
to graphs of simple or of simplicial polytopes does not make the graph
isomorphism problem easier.

\begin{theorem}
\label{thm:polyGraphIsoCompl}
The graph isomorphism problem restricted to graphs of polytopes is
graph isomorphism complete, even if one restricts the problem further
to the class of graphs of simple or to the class of graphs of
simplicial polytopes.
\end{theorem}

In particular, Theorem~\ref{thm:polyGraphIsoCompl} for simple polytopes
implies that the restricition of the graph isomorphism problem to
regular graphs is graph isomorphism complete.
As mentioned in the introduction, this is well-known. However, the
different reductions due to 
   Booth \cite{Booth78},
   Miller \cite{Miller79}, and 
   Corneil \& Kirkpatrick \cite{CorneilKirkpatrick80}
do not  produce polytopal graphs.\\

The equivalences ``(i) $\Leftrightarrow$ (iv)'' and ``(i)
$\Leftrightarrow$ (v)'' in Proposition~\ref{prop:GP} together with
part~(ii) of Remark~\ref{rem:GP} also imply the main result 
of this section.

\begin{theorem}
\label{thm:polyIsoCompl}
The polytope isomorphism problem is graph isomorphism complete, even
if one restricts the problem further to the class of simple or to the
class of simplicial polytopes.
\end{theorem}

In fact, since the duals of the polytopes $\poly{G}$ are stacked
polytopes, Theorems~\ref{thm:polyGraphIsoCompl}
and~\ref{thm:polyIsoCompl} even hold for the very restricted class of 
stacked polytopes.

\begin{figure}
  \begin{center}
    \ifusepstricks
       \begin{pspicture}(-2,-2)(2.8,3)\small

   \psset{xunit=6mm}
   \psset{yunit=7mm}

   \rput( 6, 0){\rnode{gi}{    \psframebox{$G$ isomorphic to $H$}}}   
   \rput(-2, 3){\rnode{congr}{ \psframebox{$\poly{G}$ congruent to $\poly{H}$}}} 
   \rput(-2, 1){\rnode{aff}{   \psframebox{$\poly{G}$ affinely isomorphic to $\poly{H}$}}}   
   \rput(-2,-1){\rnode{proj}{  \psframebox{$\poly{G}$ projectively isomorphic to $\poly{H}$}}}
   \rput(-2,-3){\rnode{comb}{  \psframebox{$\poly{G}$ combinatorially isomorphic to $\poly{H}$}}}

   \psset{nodesep=.5mm, doubleline=true}

   \ncline{->}{congr}{aff}
   \ncline{->}{aff}{proj}
   \ncline{->}{proj}{comb}

   \psset{nodesep=.5mm, doubleline=true, doublesep=1pt, linewidth=1.5pt}

   \nccurve[angleB=-90]{->}{comb}{gi}
   \nccurve[angleA=90]{->}{gi}{congr}

\end{pspicture}
    \else
       \epsfig{file=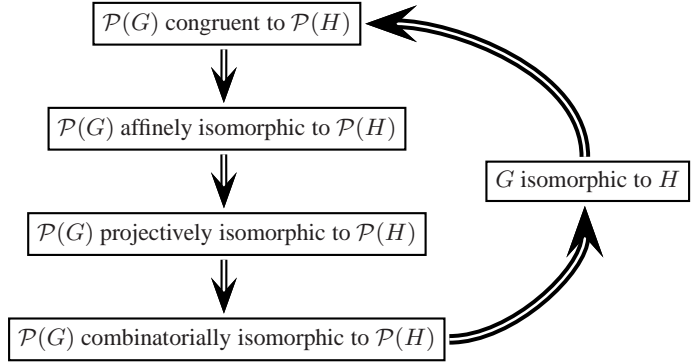,width=.8\textwidth}
    \fi
    \medskip\medskip
    \caption{Illustration of the proof of
      Theorem~\ref{thm:geomiso}. The two long arrows show implications
      following from Proposition~\ref{prop:GP} and the remarks above,
      while the three other arrows represent implications that hold in
      general.}
    \label{fig:Thm3}
  \end{center}
\end{figure}

Parts~(iii) and~(iv) of Remark~\ref{rem:GP} show that the polytope
isomorphism problem remains graph isomorphism complete if additionally
$\mathcal{V}$- and $\mathcal{H}$-descriptions of the polytopes are
provided as input data. Even more: if two graphs~$G$ and~$H$ are
isomorphic, then the polytopes~$\poly{G}$ and~$\poly{H}$ are affinely
isomorphic (here, of course, all cutting operations have to be
performed ``in the same way''). If we start our constructions with the
regular $(n-1)$-dimensional simplex 
$$
\convOp\{(1,0,\dots,0),(0,1,0,\dots,0),\dots,(0,\dots,0,1)\}\subset\R^n
$$
embedded into $\R^n$, then we can furthermore achieve that
~$\poly{G}$ and~$\poly{H}$ are congruent if~$G$ and~$H$ are
isomorphic.  Consequently, Proposition~\ref{prop:GP} (``(iv)
$\Rightarrow$ (i)'') shows that the implications depicted by the arcs
in Figure~\ref{fig:Thm3} hold.

Thus, all statements in that figure are pairwise equivalent, and
therefore, we have proved the following result on the geometric
polytope isomorphism problems.

\begin{theorem}
\label{thm:geomiso}
The graph isomorphism problem is Karp reducible to the polytope
congruence problem as well as to the affine and to the projective
polytope isomorphism problem (with respect to both $\mathcal{V}$- and
$\mathcal{H}$-descriptions). This remains true for the restrictions to
simple or to simplicial polytopes.
\end{theorem}

It has been proved by Akutsu~\cite{Akutsu98} that the graph
isomorphism problem is Karp reducible to the congruence problem for
arbitray point sets which are not necessarily in convex position.  By
an obvious modification his construction can be changed to produce
point sets in convex position, thus giving an alternative proof of the
graph isomorphism hardness of polytope congruence. The resulting
reduction also shows the graph isomorphism hardness of the affine and
of the projective polytope isomorphism problem. However, this
construction does neither yield simple nor simplicial polytopes.
Moreover, it is not clear that combinatorial isomorphism between two
constructed polytopes implies that the graphs started from are
isomorphic as well. Hence, Akutsu's construction does not provide an
alternative proof of Theorem~\ref{thm:polyIsoCompl}.

Akutsu~\cite{Akutsu98} also gave a Karp reduction of the congruence
problem to the (labeled) graph isomorphism problem, thus showing that
the former one is graph isomorphism complete. It is, however, unknown
if the affine or the projective polytope isomorphism problem can be
reduced to the graph isomorphism problem as well.



\section{Polynomiality Results for Bounded Dimension}
\label{sec:boundDim}


For polytopes of dimension one or two both the graph isomorphism
problem as well as the polytope isomorphism problem can obviously be
solved in polynomial time, even in linear time.
 It is well-known that two three-dimensional
polytopes are combinatorially isomorphic if and only if their graphs
are isomorphic (this follows from the theorem of
Whitney~\cite{Whitney35} on the uniqueness of the plane embedding of a
planar three-connected graph). Since the graph isomorphism problem for
planar graphs can be solved in linear time by an algorithm due to
Hopcroft and Wong~\cite{HopcroftWong74}, both the graph and the
polytope isomorphism problem for three-dimensional polytopes can thus
be solved in linear time.

The main result of this section is that, for every bounded dimension,
the polytope isomorphism problem can be solved in polynomial time.
First, we consider the polytope isomorphism problem restricted to
simple polytopes.

If two polytopes are combinatorially isomorphic, then they have the
same dimension as well as the same number of vertices and the same
number of facets.  It is trivial to determine these three parameters
from the vertex-facet incidences. Therefore, we will assume in the
formulations of our algorithms for the polytope isomorphism problem
that each of the three parameters has the same value for both input
polytopes.

In the following, we denote the (abstract) sets of vertices and edges
of a polytope~$P$ by $\vertices{P}$ and $\edges{P}$, respectively.


\begin{proposition}\label{prop:iso of simple d-polytopes}
  The polytope isomorphism problem for simple $d$-polytopes can be
  solved in $\order{d!\cdot d^2\cdot n^2}$ time, where~$n$ is the
  number of vertices. In particular, there is an algorithm for the
  problem whose running time is $\order{n^2}$ for bounded~$d$.
\end{proposition}

\begin{proof}
  The key observation for the proof is that every isomorphism of two simple
  polytopes is determined by its restriction to an arbitrary vertex
  and its neighborhood.
  
  Throughout the following, assume $d>2$.  Let~$P$ be a simple
  $d$-polytope. We denote by $\neighbors{x}$ the set of neighbors of
  the node~$x$ in~$\graph{P}$ and define
  $\neighborsIncl{x}:=\neighbors{x}\cup\{x\}$.  The $2$-skeleton
  of~$P$ induces, for each edge $\{v,w\}\in\edges{P}$ of~$P$, a
  bijection $ \edgebij{v}{w}\colon
  \neighborsIncl{v}\longrightarrow\neighborsIncl{w} $ with
  $\edgebij{v}{w}(v)=v$, $\edgebij{v}{w}(w)=w$, and
  $\edgebij{v}{w}(u)$ being the other (than~$v$) neighbor of~$w$ in
  the $2$-face spanned by $v$, $w$, and~$u$.
  
  Suppose $\pi\colon\vertices{P}\longrightarrow\vertices{Q}$ is an
  isomorphism of the graphs~$\graph{P}$ and~$\graph{Q}$ of two simple
  $d$-polytopes~$P$ and~$Q$, respectively. For each $v\in\vertices{P}$
  we denote by $ \pi_v\colon\ 
  \neighborsIncl{v}\longrightarrow\vertices{Q} $ the restriction
  of~$\pi$ to $\neighborsIncl{v}$. For every node $v\in\vertices{P}$,
  we have
  \begin{equation}
    \label{eq:i}
    \pi_v(\neighbors{v})=\neighbors{\pi_v(v)}\enspace,
  \end{equation}
  and, for every edge $\{v,w\}\in\edges{P}$,
  \begin{equation}
    \label{eq:ii}
    \pi_w=\edgebij{\pi_v(v)}{\pi_v(w)} \circ \pi_v \circ \edgebij{w}{v}
  \end{equation}
  (where the composition $\edgebij{\pi_v(v)}{\pi_v(w)} \circ
  \pi_v$ is well-defined due to~(\ref{eq:i})).
  
  Conversely, for two simple $d$-polytopes~$P$ and~$Q$ with~$n$
  vertices consider any set of maps $\pi_v\colon
  \neighborsIncl{v}\longrightarrow\vertices{Q}$ ($v\in\vertices{P}$).
  If the maps are consistent, i.e.,
  \begin{equation}
    \label{eq:iii}
    \pi_v(u)=\pi_w(u)\qquad
    (v,w,u\in\vertices{P},
    u\in\neighborsIncl{v}\cap\neighborsIncl{w})\enspace,
  \end{equation}
  then there is a unique map $\pi\colon
  \vertices{P}\longrightarrow\vertices{Q}$ such that $\pi_v$ is the
  restriction of~$\pi$ to $\neighborsIncl{v}$ for all
  $v\in\vertices{P}$. We claim that, if $\pi$ satisfies~(\ref{eq:i})
  for all $v\in\vertices{P}$, then~$\pi$ is an isomorphism of
  $\graph{P}$ and $\graph{Q}$ (and thus a combinatorial isomorphism
  between~$P$ and~$Q$). To see this, it suffices to show that~$\pi$ is
  surjective, since both~$\graph{P}$ and~$\graph{Q}$ have the same
  number~$n$ of nodes. Suppose that~$\pi$ is not surjective.
  Since~$\graph{Q}$ is connected, there is a node~$y$ of~$\graph{Q}$
  which is not contained in the image of~$\pi$ and which has a
  neighbor~$x$ that is the image~$\pi(v)$ of some node~$v$
  of~$\graph{P}$. However, this contradicts~(\ref{eq:i}).
  
  Thus, we may check two simple $d$-polytopes~$P$ and~$Q$ both
  with~$n$ vertices, given by their vertex-facet incidences, for
  combinatorial isomorphism in the following way. First, we compute
  the graphs~$\graph{P}$ and~$\graph{Q}$, as well as the bijections
  $\edgebij{v}{w}$ ($\{v,w\}\in\edges{P}$) and $\edgebij{x}{y}$
  ($\{x,y\}\in\edges{Q}$).  Furthermore, a node $v_0\in\vertices{P}$
  is fixed together with a spanning tree~$T_0$ of~$\graph{P}$, rooted
  at~$v_0$.
  
  Then, for each node $x\in\vertices{Q}$ and for each bijection
  $\pi_{v_0}\colon \neighborsIncl{v_0}\longrightarrow\neighborsIncl{x}$
  with $\pi_{v_0}(v_0)=x$, we perform the following steps:
  \begin{enumerate}
  \item Compute~$\pi_w$ ($w\in\vertices{P}$) ``along~$T_0$'' by means
    of~(\ref{eq:ii}); if, for some~$w$, condition~(\ref{eq:i}) is not
    satisfied, then continue with the next bijection~$\pi_{v_0}$.
  \item If~(\ref{eq:iii}) is not satisfied, then continue with the
    next bijection~$\pi_{v_0}$.
  \item Construct the map~$\pi$ (as above); STOP.
  \end{enumerate}
  
  Note that, when $\pi_w$ is computed in Step~1 from the parent~$v$
  of~$w$ in~$T_0$, condition~(\ref{eq:i}) is satisfied for~$v$ (thus,
  the composition in~(\ref{eq:ii}) is well-defined). It follows from
  the discussion above that the two polytopes are isomorphic if and
  only if the algorithm stops in Step~3; in this case, the constructed
  $\pi$ is an isomorphism between~$P$ and~$Q$. 
  
  Computing the graphs~$\graph{P}$ and~$\graph{Q}$, as well as the
  bijections $\edgebij{v}{w}$ ($\{v,w\}\in\edges{P}$) and
  $\edgebij{x}{y}$ ($\{x,y\}\in\edges{Q}$), can be performed in
  $\order{d \cdot n^2}$ time, if we perform the
  following preprocessing in advance. From the vertex-facet incidences
  compute, for each vertex, a sorted list of indices of the facets
  containing this vertex.  This can be performed in $\order{n^2}$
  steps (note that simple polytopes never have more facets than
  vertices).  Compute a similar incidence list for each facet.
  
  Moreover, using this data structure, none of the four steps in the
  for-loop needs more than $\order{d^2\cdot n}$ time (the critical
  part being Step~2).  Since the body of the for-loop is not executed
  more than $n\cdot d!$ times, this yields an $\order{d!\cdot d^2\cdot
    n^2}$ time algorithm.\hspace*{\fill}\qed
  \end{proof}

Luks~\cite{Luks82} gave a polynomial time algorithm for the graph
isomorphism problem on graphs of bounded maximal degree. Since the
graph of a simple $d$-polytope is $d$-regular, his algorithm runs in
polynomial time on graphs of simple polytopes of bounded dimension.

\begin{proposition}
\label{prop:graphIsoSimpleFix}
The isomorphism problem for graphs of simple polytopes of bounded
dimension can be solved in polynomial time.
\end{proposition}

As two simple polytopes are isomorphic if and only if their graphs are
isomorphic (cf. the proof of Proposition~\ref{prop:GP}) and since it
is easy to compute efficiently the graph of a polytope from its
vertex-facet incidences, Proposition~\ref{prop:graphIsoSimpleFix}
implies the polynomiality statement of Proposition~\ref{prop:iso of
  simple d-polytopes}. The algorithm described in the proof of
Proposition~\ref{prop:iso of simple d-polytopes}, however, is both
much simpler to understand/implement, and much faster than Luks'
method. The running time of the original version of Luks' isomorphism
test for graphs of bounded maximal degree~$d$ is~$n^{\order{d^3}}$,
where~$n$ is the number of nodes. According to Luks, this could be
improved to $n^{O(d/\log d)}$~\cite{Luks01}, which, in general,
is much larger than $\order{n^2}$.

Two polytopes are isomorphic if and only if their dual polytopes are
isomorphic. Since the transpose of a vertex-facet incidence matrix of
a polytope is a vertex-facet incidence matrix of the dual polytope,
Proposition~\ref{prop:iso of simple d-polytopes} implies its own
analogue for simplicial polytopes.

\begin{proposition}\label{prop:iso of simplicial d-polytopes}
  The polytope isomorphism problem for simplicial $d$-polytopes can be
  solved in $\order{d!\cdot d^2\cdot m^2}$ time, where~$m$ is the
  number of facets. In particular, there is an algorithm for the
  problem whose running time is $\order{m^2}$ for bounded~$d$.
\end{proposition}

While simple and simplicial polytopes play symmetric roles with
respect to the polytope isomorphism problem, they may play different
roles with respect to the graph isomorphism problem. In particular, it
is unknown whether Proposition~\ref{prop:graphIsoSimpleFix} is also
true for the graphs of simplicial polytopes of bounded dimensions (see
Sect.~\ref{sec:conclu}).

The main result of this section is that one can extend
Propositions~\ref{prop:iso of simple d-polytopes} and~\ref{prop:iso of
  simplicial d-polytopes} to arbitrary polytopes (of bounded
dimension).

\begin{theorem}\label{thm:iso of d-polytopes}
  The polytope isomorphism problem can be solved in
  \begin{center}
    $\order{d^2\cdot 2^{d^2}\cdot\min\{n,m\}^{d^2/2}}$\hspace{.5cm}time, 
  \end{center}
  where $d$, $n$ and~$m$ are
  the dimension, the number of vertices, and the number of facets,
  respectively.  In particular, there is an algorithm for the problem
  whose running time is $\order{\min\{n,m\}^{d^2/2}}$ for bounded~$d$.
\end{theorem}

\begin{proof}
  The core of the proof is a result of Bayer~\cite[see Theorem~3 and
  its proof]{Bayer88} implying that every combinatorial isomorphism
  between the barycentric subdivisions of two polytopes~$P$ and~$Q$
  induces a combinatorial isomorphism between~$P$ and~$Q$ or its dual
  polytope~$\dualPolytope{Q}$, where the \emph{barycentric
    subdivision} of a $d$-polytope~$R$ may be defined as any
  $d$-polytope~$\bary{R}$, whose vertices correspond to the
  non-trivial faces ($\not=\varnothing, R$) of~$R$, and whose facets
  correspond to the maximal chains in the face lattice of~$R$. Such a
  polytope~$\bary{R}$ can be constructed from~$R$ by performing
  stellar subdivisions of all $k$-faces for $k=d-1,\dots,1$ (see
  Fig.~\ref{fig:bary}).
 
\begin{figure}
  \centering
  \subfigure[$R$]{
    \epsfig{figure=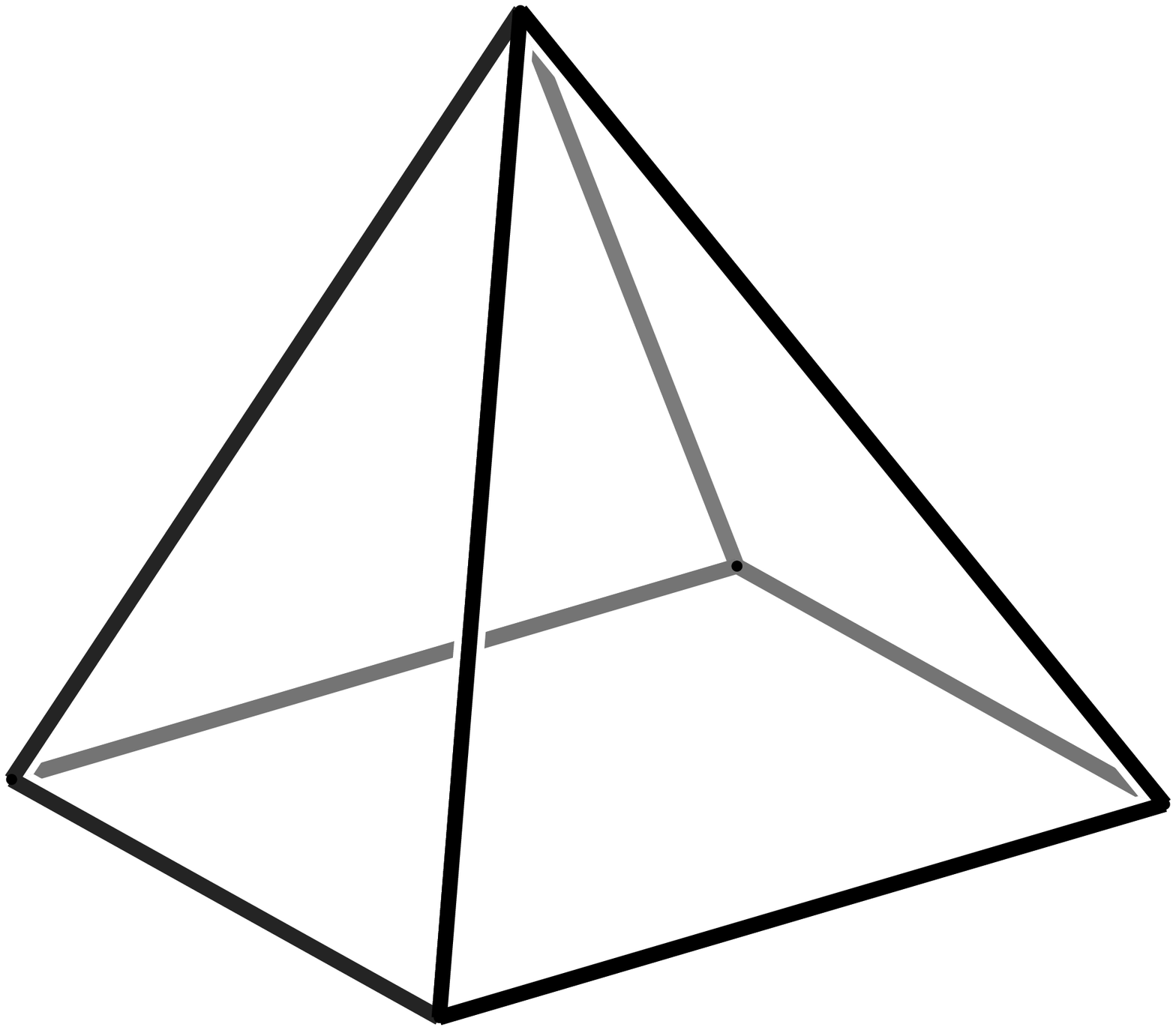,width=.3\textwidth}} \hfill
  \subfigure[$\bary{R}$]{
    \epsfig{figure=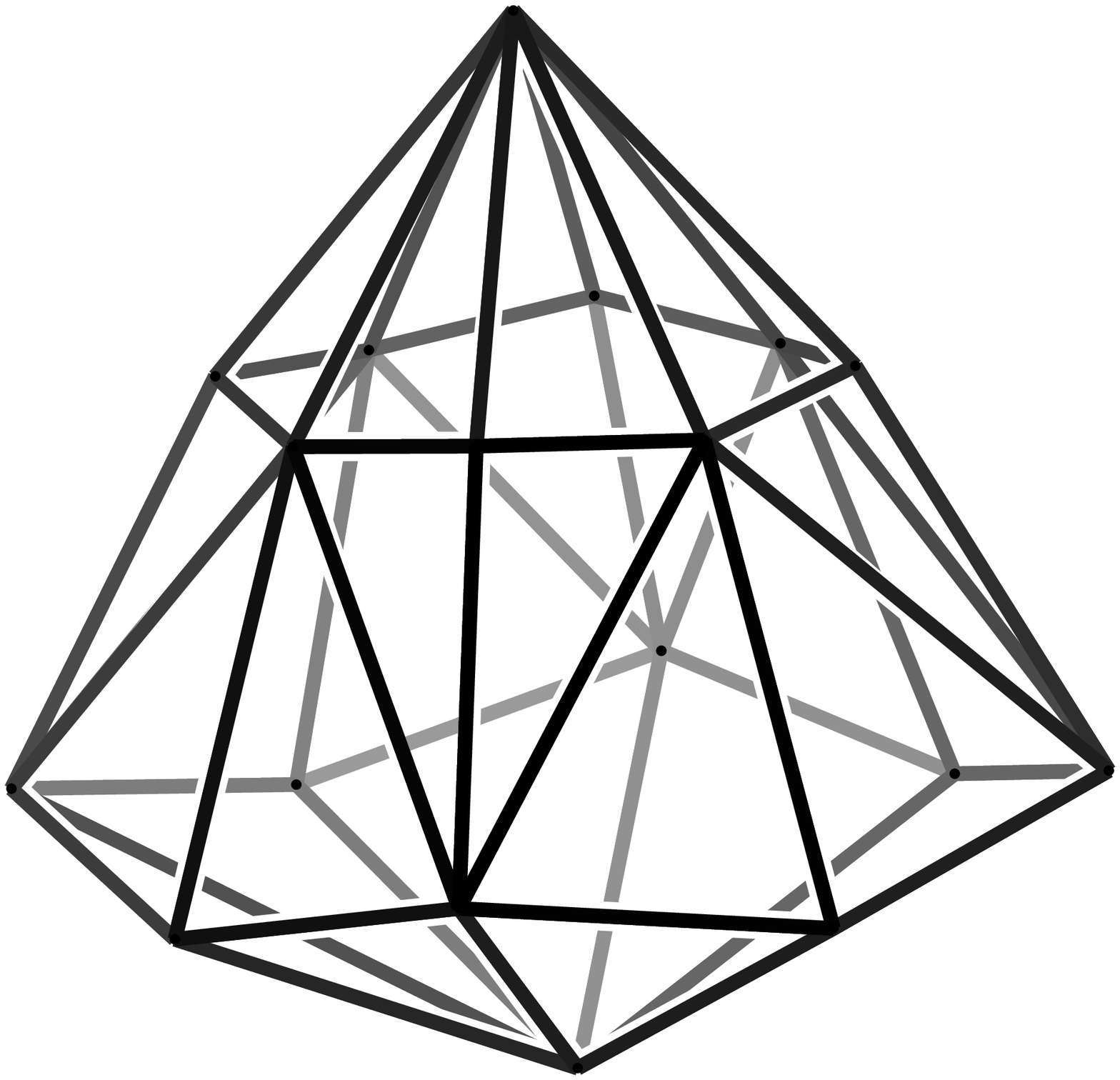,width=.3\textwidth}} \hfill
  \subfigure[$\dualPolytope{\bary{R}}$]{
    \epsfig{figure=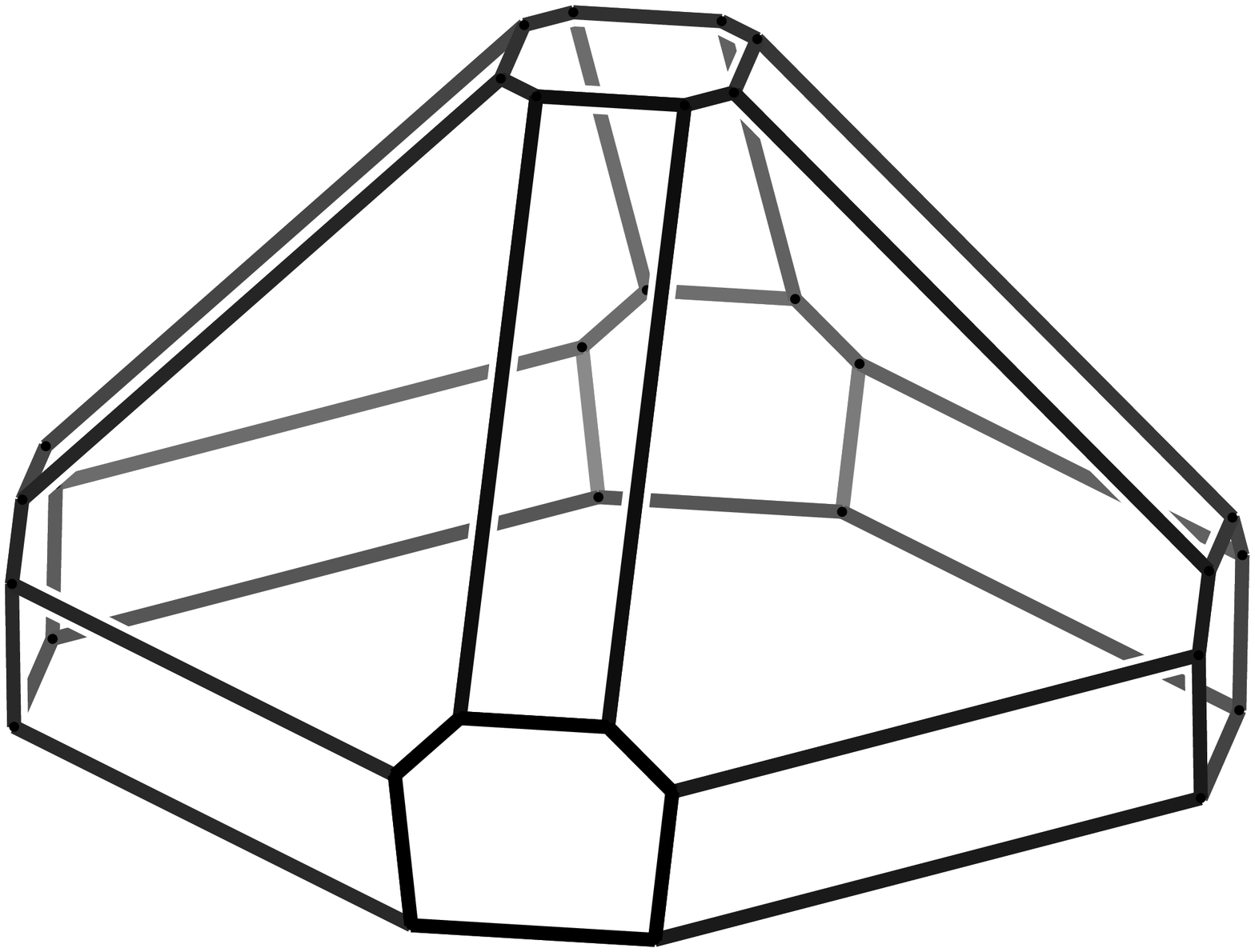,width=.3\textwidth}}
  \caption{The barycentric subdivision $\bary{R}$ of a 3-dimensional
    pyramid~$R$ and  the dual polytope $\dualPolytope{\bary{R}}$.}
  \label{fig:bary}
\end{figure}
 
  The dual~$\dualPolytope{\bary{R}}$ is a simple $d$-polytope whose
  vertices correspond to the maximal chains in the face lattice
  of~$R$, and whose facets correspond to the (non-empty) faces of~$R$.
  In the graph of~$\dualPolytope{\bary{R}}$, two nodes are adjacent if
  and only if the corresponding two maximal chains of~$R$ arise from
  each other by exchanging the two $i$-dimensional faces (for some
  unique~$0\leq i\leq d-1$). This defines a labeling of the edges
  of~$\graph{\dualPolytope{\bary{R}}}$ by $\{0,1,\dots,d-1\}$, where
  for each node the~$d$ incident edges receive pairwise different
  labels.  Thus, Bayer's result implies that two polytopes~$P$ and~$Q$
  are combinatorially isomorphic if and only if there is a
  combinatorial isomorphism between~$\dualPolytope{\bary{P}}$
  and~$\dualPolytope{\bary{Q}}$ that induces a label preserving
  isomorphism of $\graph{\dualPolytope{\bary{P}}}$ and
  $\graph{\dualPolytope{\bary{Q}}}$.
  
  Let~$P$ and~$Q$ be two $d$-polytopes, given by their vertex-facet
  incidences, which are to be checked for combinatorial isomorphism.
  We first compute the Hasse diagram of the face lattice of each of
  the two polytopes (with nodes labeled by the dimensions of the
  respective faces). From the Hasse diagrams we then enumerate the
  maximal chains of the face lattices of~$P$ and~$Q$, the vertex-facet
  incidences of~$\dualPolytope{\bary{P}}$
  and~$\dualPolytope{\bary{Q}}$, as well as the
  graphs~$\graph{\dualPolytope{\bary{P}}}$ and
  $\graph{\dualPolytope{\bary{Q}}}$ together with their edge
  labelings.
  
  Once these data are available, we can
  check~$\dualPolytope{\bary{P}}$ and~$\dualPolytope{\bary{Q}}$ for
  combinatorial isomorphism as in the proof of
  Proposition~\ref{prop:iso of simple d-polytopes}. However, since we
  only allow isomorphisms that respect the edge labeling
  of~$\graph{\dualPolytope{\bary{P}}}$ and
  $\graph{\dualPolytope{\bary{Q}}}$, for each potential image~$x$ of
  the fixed vertex~$v_0$ of~$\dualPolytope{\bary{P}}$ there is only
  one possible bijection between
  $\neighborsIncl{v_0}\longrightarrow\neighborsIncl{x}$ (with the
  notation adapted from the proof of Proposition~\ref{prop:iso of
    simple d-polytopes}). Consequently, the running time of
  checking~$\dualPolytope{\bary{P}}$ and~$\dualPolytope{\bary{Q}}$ for
  (label preserving) combinatorial isomorphism can be estimated by
  $\order{d^2\cdot \zeta^2}$, where~$\zeta$ is the number of maximal
  chains in the face lattice of~$P$ or~$Q$. We may assume that~$\zeta$
  is the same number for both~$P$ and~$Q$, since otherwise we know
  already that~$P$ is not isomorphic to~$Q$; the same holds for the
  numbers~$\varphi$, $\alpha$, $n$, and~$m$ arising below.
  
  Computing the Hasse diagrams of~$P$ and~$Q$ can be performed in time
  $\order{\varphi\cdot\alpha\cdot\min\{n,m\}}$ (see
  \cite{KaibelPfetsch02}), where $\varphi$, $\alpha$, $n$, and~$m$ are
  the number of faces, the number of vertex-facet incidences, the
  number of vertices, and the number of facets of~$P$ or~$Q$,
  respectively.  Enumerating all maximal chains in the face lattices
  of~$P$ and~$Q$, can be done in time proportional to $d\cdot\zeta$.
  Computing the vertex-facet incidences of~$\dualPolytope{\bary{P}}$
  and~$\dualPolytope{\bary{Q}}$ takes no more than
  $\order{\varphi\cdot\zeta}$ steps. The edge labeled
  graphs~$\graph{\dualPolytope{\bary{P}}}$ and
  $\graph{\dualPolytope{\bary{Q}}}$ as well as the maps
  $\Psi_{\star,\star}$ can be obtained (see the proof of
  Proposition~\ref{prop:iso of simple d-polytopes}) in time
  $\order{d\cdot\zeta^2}$ (note that $\Psi_{\star,\star}$ can easily
  be obtained from the edge labeling).
  
  By the upper bound theorem for convex polytopes (see, e.g.,
  \cite[Thm.~8.23]{Ziegler95}) applied to the barycentric
  subdivisions, we have $\zeta\in\order{\varphi^{ d/2}}$, where, again
  by the upper bound theorem, we can estimate $\varphi\le
  2^d\cdot\min\{n,m\}^{d/2}$. In total, the running time of the sketched
  algorithm thus can be bounded by
  $\order{d^2\cdot\zeta^2}=\order{d^2\cdot
    2^{d^2}\cdot\min\{n,m\}^{d^2/2}}$.
 \hspace*{\fill}\qed
\end{proof}

The estimate for the running time in this proof uses the upper bound
theorem twice in a nested way. Thus, in practice the running time will
usually be much smaller than the estimate of the theorem.\\

All three geometric polytope isomorphism problems are polynomially
solvable in bounded dimensions: If the vertices
$\coordVertices{P}\subset\R^d$ and $\coordVertices{Q}\subset\R^d$ of two
polytopes~$P$ and~$Q$ are given, one may first choose a maximal
affinely independent set~$S\subset\coordVertices{P}$. Then, for each
maximal affinely independent set $T\subset\coordVertices{Q}$ and for each
bijection~$\pi\colon S\longrightarrow T$ one checks whether the affine
map $\aff{P}\longrightarrow\aff{Q}$ induced by~$\pi$ satisfies
$\pi(\coordVertices{P})=\coordVertices{Q}$.  Similar algorithms can be
constructed for checking congruence or projective isomorphism. If the
polytopes are specified by $\mathcal{H}$-descriptions one may first
compute $\mathcal{V}$-representations (in polynomial time, since the
dimension is bounded).

There are more elaborate algorithms for checking two sets of~$n$
points in~$\R^d$ for congruence. For instance, Brass \&\ 
Knauer~\cite{BrassKnauer02} describe a deterministic $\order{n^{\lceil
    d/3\rceil}\cdot\log n}$ algorithm.




\section{Conclusions}
\label{sec:conclu}

\begin{figure}
    \begin{center}  
      \ifusepstricks
         \newcommand{\problembox}[1]{\psframebox[fillstyle=solid,fillcolor=white]{\tabcolsep=0pt\begin{tabular}[t]{c}#1\end{tabular}}}
\newcommand{\lproblembox}[1]{\psframebox[fillstyle=solid,fillcolor=white]{\parbox{2.7cm}{\centerline{\tabcolsep=0pt\begin{tabular}[t]{c}#1\end{tabular}}}}}
\newcommand{\sproblembox}[1]{\psframebox[fillstyle=solid,fillcolor=white]{\parbox{1.5cm}{\tabcolsep=0pt\begin{tabular}[t]{c}#1\end{tabular}}}}

\begin{pspicture}(-1.24,-3.5)(4.4,2.5)\small

   \psset{xunit=31mm}
   \psset{yunit=8mm}

   \newrgbcolor{yellowl}{1 1 .5}

   \pspolygon[linewidth=.1pt,linearc=0.3,fillstyle=solid,fillcolor=yellowl]%
      (-.58,-4)(2.2,-4)(2.2,4)(-.58,4)(-.58,2)(.6,2)(.6,-2)(-.58,-2)
   \put(2.5,-3){\emph{\textbf{graph isomorphism complete}}}

   \rput( 1.02,   0){\rnode{gi}{     \problembox{graph\\iso problem}}}   
   \rput( 1.8, 0){\rnode{lgi}{    \problembox{labeled graph\\iso problem}}}
 
   \rput( 0,   3){\rnode{congr}{  \lproblembox{polytope congruence\\problem}}}    
   \rput( 0,   1){\rnode{aff}{    \lproblembox{affine polytope\\iso problem}}}    
   \rput( 0,  -1){\rnode{proj}{   \lproblembox{projective polytope\\iso problem}}}    
   \rput( 0,  -3){\rnode{pi}{     \lproblembox{(comb.) polytope\\iso problem}}}    
   \rput(-1,  -2){\rnode{fli}{    \sproblembox{face-lattice\\iso problem}}}    
   \rput(-1,  -4){\rnode{sd}{     \sproblembox{self-duality\\problem}}}    

   \psset{nodesep=.7mm}
  
   \ncarc{->}{gi}{lgi}
   \ncarc{->}{lgi}{gi}

   \ncarc[angleA=-140]{<-}{lgi}{pi}
   \ncline{->}{fli}{pi}
   \ncline{->}{sd}{pi}

   \psset{linewidth=1.5pt}
   \ncline[angleA=-140]{->}{gi}{aff}
   \ncline[angleA=-140]{->}{gi}{proj}
   \ncarc[angleA=-140]{->}{gi}{pi}

   \psset{linewidth=1.5pt,linestyle=dashed}
   \ncarc[angleB=90]{->}{congr}{lgi}
   \ncarc[angleA=140]{<-}{congr}{gi}
\end{pspicture}
      \else
         \epsfig{file=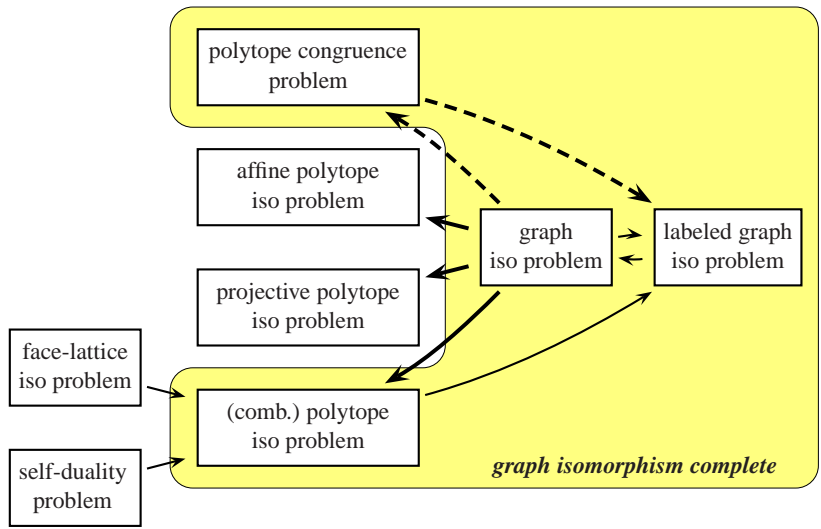,width=\textwidth}
      \fi
      \caption{Illustration of several decision problems treated in
        the paper; each arrow means that there is a Karp-reduction
        from the problem at its tail to the problem at its head. The
        solid bold arrows indicate new reductions and the two dashed ones
        emerge from Akutsu~\cite{Akutsu98}.}  
      \label{fig:BigPicture}
    \end{center} 
\end{figure}

The main results of this paper are on the one hand the graph
isomorphism completeness of the general (combinatorial) polytope
isomorphism problem and, on the other hand, the fact that this problem
can be solved in polynomial time if the dimensions of the polytopes
are bounded by a constant (see Table~\ref{tab:summary} for an overview
of the complexity results and Fig.~\ref{fig:BigPicture} for a sketch
of the complexity theoretic landscape considered in this paper).

\begin{table}
  \caption{Summary of the complexity results (GI means graph
    isomorphism).}
  \label{tab:summary}
    \newcommand{\rb}[1]{\raisebox{1ex}[-1ex]{#1}}

      

\begin{tabular*}{\linewidth}{@{\extracolsep{\fill}}llccc@{}}\toprule    
                    &            & \begin{tabular}{c}
                                       polytopal\\graph isomorphism
                                    \end{tabular} 
                                 & \multicolumn{2}{c}{polytope isomorphism}\\
   \cmidrule{4-5}
                    &            &                   & incidences    &face lattice\\ 
   \midrule     
                    & simplicial & \emph{open}      &  polynomial   &  polynomial \\
     \rb{bounded}   & simple     &  polynomial       &  polynomial   &  polynomial \\
     \rb{dimension} & arbitrary  &  \emph{open}     &  polynomial   &  polynomial \\
   \addlinespace
                    & simplicial &     GI~complete & GI~complete & \emph{open} \\
     \rb{arbitrary} & simple     &     GI~complete & GI~complete & \emph{open} \\
     \rb{dimension} & arbitrary  &     GI~complete & GI~complete &
    \emph{open} \\ 
  \bottomrule  
\end{tabular*}
\end{table}

Our hardness result for arbitrary dimensions leaves open the question
for the complexity status of the problem to decide whether two
polytopes given by their entire face lattices (rather than by their
vertex-facet incidences only) are isomorphic (\emph{face lattice
  isomorphism problem}).  The algorithm described in the proof of
Theorem~\ref{thm:iso of d-polytopes} does not solve the face lattice
isomorphism problem in polynomial time, since the number of maximal
chains in the face lattice of a polytope is not bounded polynomially
in the size of the face lattice.

The polynomiality result on the polytope isomorphism problem in
bounded dimensions (Theorem~\ref{thm:iso of d-polytopes})
theoretically supports the empirical evidence that checking whether
two polytopes of moderate dimensions are isomorphic is not too hard.
This evidence stems from applying McKay's \nauty\ to the vertex-facet
incidences as mentioned in the introduction. It may be that one can
turn our algorithm into a computer code that becomes compatible with
\nauty\ for checking combinatorial polytope isomorphism.

The remaining two \emph{open} entries in Table~\ref{tab:summary}
concern the complexity of the graph isomorphism problem restricted to
graphs of arbitrary (or simplicial) polytopes of bounded dimensions. A
polynomial time algorithm for this problem would perhaps not be as
interesting as the potential result that the problem is graph
isomorphism complete, because the latter result would show that the
class of graphs of polytopes of any fixed dimension is in a sense
``structurally as rich'' as the class of all graphs.

If the face lattice isomorphism problem turned out to be graph
isomorphism complete, then this would immediately imply that the face
lattice isomorphism problem is also ``(finite) poset isomorphism
complete.''  Similarly to the case of the graphs of polytopes of fixed
dimensions, one might interpret such a result as an evidence that the
class of face lattices form a ``structurally rich'' subclass of the
class of all (finite) posets. 

An interesting special case of the combinatorial polytope isomorphism
problem is the \emph{polytope self-duality problem}, which asks to
decide whether a polytope, given by its vertex-facet incidences, is
combinatorially isomorphic to its dual. Clearly, by
Theorem~\ref{thm:iso of d-polytopes} this problem can be solved in
polynomial time for bounded dimensions. However, its general
complexity status remains open.


\begin{acknowledgement}
  We are indebted to G\"unter M.~Ziegler for careful reading an
  earlier version of the paper and for stimulating discussions, in
  particular, for bringing Bayer's work~\cite{Bayer88} to our
  attention.  We thank G\"unter Rote for inspiring comments, Christian
  Knauer for pointing to Akutsu's results \cite{Akutsu98}, and
  Christoph Eyrich for valuable \LaTeX\ advice.
\end{acknowledgement}






\begin{thebibliography}{10}

\bibitem{Akutsu98}
Akutsu, T.: On determining the congruence of point sets in {$d$} dimensions.
\newblock Comput.~Geom.~{\bf 9}, 247--256 (1998)

\bibitem{Babai95}
Babai, L.: Automorphism groups, isomorphism, reconstruction. 
In: R.L.~Graham et al.~(eds.): Handbook of Cominatorics (pp.~1447--1540), 
    Amsterdam, North-Holland: Elsevier 1995

\bibitem{Bayer88}
Bayer, M.: Barycentric subdivisions.
\newblock Pac. J. Math.~{\bf 135}, 1--17 (1988)

\bibitem{BlindMani87}
Blind, R. and Mani-Levitska, P.: On puzzles and polytope isomorphisms.
\newblock Aequationes Math.~{\bf 34}, 287--297 (1987)

\bibitem{Booth78}
Booth, K.: Isomorphism testing for graphs, semigroups, and finite automata are
  polynomially equivalent problems.
\newblock SIAM J.~Comput.~{\bf 7}, 273--279 (1978)

\bibitem{BoothLueker75}
Booth, K. and Lueker, G.~S.: Linear algorithms to recognize interval graphs and
  test for the consecutive ones property.
\newblock Proc.~7th Ann.~ACM Symp.~Theory Comput., 255--265 (1975)

\bibitem{BrassKnauer02}
Brass, P. and Knauer, C.: Testing the congruence of $d$-dimensional point sets.
\newblock J.~Comput.~Geom.~Appl.~{\bf 12}, 115--124 (2002)

\bibitem{ColbournColbourn78}
Colbourn, M.J. and Colbourn, J.: Graph isomorphism and self-complementary graphs.
\newblock SIGACT News~{\bf 10}, 25--30 (1978)

\bibitem{CorneilKirkpatrick80}
Corneil, D. and Kirkpatrick, D.: A theoretical analysis of various heuristics
  for the graph isomorphism problem.
\newblock SIAM J.~Comput.~{\bf 9}, 281--297 (1980)

\bibitem{Fortin96}
Fortin, S.: The graph isomorphism problem, tech.~report TR 96-20, Dept.~of
  Computer Science, University of Alberta, Edmonton, Alberta, Canada, 1996
\newblock
  \\\url{ftp://ftp.cs.ualberta.ca/pub/TechReports/1996/TR96-20/TR96-20.ps.gz}

\bibitem{GawrilowJoswig00b}
Gawrilow, E. and Joswig, M.: \newblock polymake: {A} framework for analyzing convex
  polytopes. 
\newblock In:  G.~Kalai and G.M.~Ziegler (eds.): 
\newblock      Polytopes --- Combinatorics and Computation (DMV-Seminars, pp.~43--74)
\newblock      Basel: Birkh\"auser-Verlag Basel 2000,\\
\newblock see also \url{http://www.math.tu-berlin.de/diskregeom/polymake}

\bibitem{Gruenbaum67}
Gr{\"u}nbaum, B.: Convex Polytopes. New York: Interscience Publishers John Wiley \& Sons 1967


\bibitem{Harary69}
Harary, F.: Graph Theory.
\newblock Reading, MA: Addison-Wesley 1969

\bibitem{Hoffmann82}
Hoffmann, C.: Group-theoretic Algorithms and Graph Isomorphism (Lecture Notes
  in Computer Science~136) Berlin-Heidelberg-New York: Springer 1982

\bibitem{HopcroftWong74}
Hopcroft, J. and Wong, J.: A linear time algorithm for isomorphism of planar
  graphs.
\newblock Proc.~6th Ann.~ACM Symp.~Theory Comput., 172--184 (1974)

\bibitem{KaibelPfetsch02}
Kaibel, V. and Pfetsch, M.E.: Computing the face lattice of a polytope from
  its vertex-facet incidences.
\newblock Comput.~Geom.~(to appear)

\bibitem{Kalai88}
Kalai, G.: A simple way to tell a simple polytope from its graph.
\newblock J.~Comb.~Theory Ser.~A~{\bf 49}, 381--383 (1988)

\bibitem{KoeblerSchoeningToran93}
K{\"o}bler, J., Sch{\"o}ning, U., and Tor{\'a}n, J.: The Graph Isomorphism
  Problem: Its Structural Complexity (Progress in Theoretical Computer Science)
  Boston, MA: Birkh\"auser 1993

\bibitem{Kozen77}
Kozen, D.: On the complexity of finitely presented algebras.
\newblock Proc.~9th Ann.~ACM Symp.~Theory Comput., 164--177 (1977)

\bibitem{Kutz02}
Kutz, M.: Computing roots of directed graphs is graph isomorphism hard, tech.
  report, FU Berlin, 2002.
\newblock \texttt{arXiv:math.CO/0207020}

\bibitem{Lubiw81}
Lubiw, A.: Some {$\NP$}-complete problems similar to graph isomorphism.
\newblock SIAM J.~Comput.~{\bf 10}, 11--21 (1981)

\bibitem{Luks82}
Luks, E.: Isomorphism of graphs of bounded valence can be tested in polynomial
  time.
\newblock J.~Comput.~System~Sci.~{\bf 25}, 42--65 (1982)

\bibitem{Luks01}
Luks, E.: Personal communication, May 2001

\bibitem{Mathon79}
Mathon, R.: A note on the graph isomorphism counting problem.
\newblock Inf.~Process.~Lett.~{\bf 8}, 131--132 (1979)

\bibitem{McKay81}
McKay, B.: Practical graph isomorphism.
\newblock Numerical mathematics and computing, Proc.~10th Manitoba Conf.,
  Winnipeg/Manitoba 1980, Congr.~Numerantium 30, 45--87 (1981),
\newblock see also \url{http://cs.anu.edu.au/people/bdm/nauty}

\bibitem{Miller79}
Miller, G.L.: Graph isomorphism, general remarks.
\newblock J.~Comput.~System~Sci.~{\bf 18}, 128--142 (1979)

\bibitem{PolthierKhademPreussReitebuch01}
Polthier, K., Khadem, S., Peu{\ss}, E. and Reitebuch, U.:
\newblock Publication of interactive visualizations with \texttt{JavaView}.
\newblock To appear in: J.~Borwein et al.~(eds.): 
\newblock Multimedia Tools for Communicating Mathematics, Springer 2001,\\
\newblock see also \url{http://www.javaview.de/}

\bibitem{ReadCorneil77}
Read, R. and Corneil, D.: The graph isomorphism disease.
\newblock J.~Graph Theory~{\bf 1}, 339--363 (1977)

\bibitem{Shawe-TaylorPisanski94}
Shawe-Taylor, J. and Pisanski, T.: Homeomorphism of $2$-complexes is graph
  isomorphism complete.
\newblock SIAM J.~Comput.~{\bf 23}, 120--132 (1994)

\bibitem{Whitney35}
Whitney, H.: Congruent graphs and the connectivity of graphs.
\newblock Am.~J.~Math.~{\bf 54}, 150--168 (1932)

\bibitem{ZalcsteinFranklin86}
Zalcstein, Y. and Franklin, S.P.: Testing homotopy equivalence is isomorphism
  complete.
\newblock Discrete Appl.~Math.~{\bf 13}, 101--104 (1986)

\bibitem{Ziegler95}
Ziegler, G.M.: Lectures on Polytopes (Graduate Texts in Mathematics 152) New York:
  Springer 1995, revised edition 1998

\end{thebibliography}
\end{document}